\documentclass[12pt,a4paper]{article}
\usepackage{amsmath}
\usepackage{amssymb}
\newtheorem{theorem}{Theorem}
\newtheorem{corollary}[theorem]{Corollary}
\newtheorem{definition}[theorem]{Definition}
\newtheorem{lemma}[theorem]{Lemma}
\newtheorem{proposition}[theorem]{Proposition}
\newtheorem{remark}[theorem]{Remark}
\newenvironment{proof}[1][Proof]{\noindent\textbf{#1.} }{\ \rule{0.5em}{0.5em}}

\begin{document}

\title{A NON-STATIONARY PROBLEM COUPLING PDEs AND ODEs
MODELIZING AN AUTOMOTIVE CATALYTIC CONVERTER}

\author{JEAN-DAVID HOERNEL}
\date{}

\maketitle

\begin{center}
Universit\'e de Haute Alsace, Laboratoire de Math\'ematiques, Informatique et Applications
\vskip 0,1cm
4 rue des fr\`eres Lumi\`ere, 68093 MULHOUSE (France)
\vskip 0,1cm
j-d.hoernel@wanadoo.fr
\end{center}

\footnotetext{AMS 2000 Subject Classification: 35M10, 35K65, 35Q80}
\footnotetext{Key words and phrases: Catalytic combustion, Parabolic equation, Ordinary differential equation, Coupled system, Transmission problem.}

\begin{abstract} In this paper we prove the existence and uniqueness of the solution of a
non-stationary problem that modelizes the behaviour of the concentrations and the temperature of gases going through a cylindrical passage of an automotive catalytic converter.
This problem couples parabolic partial differential equations in a domain 
with one parabolic partial differential equation and some ordinary differential
equations on a part of its boundary.
\end{abstract}

\section{Introduction}

The starting point of this work is a non-stationary model of catalytic
converter with cylindrical passages due to Ryan, Becke and Zygourakis
[12] that is an extension to the one established by Oh and
Cavendish [10,11]. It describes the spatial and temporal evolutions of the concentrations
of $N-1$ different chemical species in gazeous phase going through a cylinder and that of the temperature in the cylinder and on its boundary. A stationary problem close to our has been
studied by Carasso [2]. Because of its internal symmetry the cylinder may be reduced to the plane domain $\Omega=\left]  0,1\right[  \times\left]  0,1\right[  $ the boundary of which
is $\Sigma=\left\{  1\right\}  \times\left]  0,1\right[  $. For $i\in\{1,...,N-1\}$, the concentrations
(resp. the temperature) inside the cylinder $\Omega$ are named $C_{if}$
(resp.$\;C_{Nf}$) and the concentrations (resp. the temperature) on
the boundary $\Sigma$ are named $C_{is}$ (resp. $C_{Ns}$). The problem is
written in a normalized way as
\begin{equation}
\left\{
\begin{array}
[c]{rcl}%
\dfrac{\partial C_{if}}{\partial z}\left(  r,z,t\right)   & = & \dfrac
{\beta_{if}}{r\left(  1-r^{2}\right)  }\dfrac{\partial}{\partial r}\left(
r\dfrac{\partial C_{if}}{\partial r}\right)  \left(  r,z,t\right)  ,\\
\dfrac{\partial C_{is}}{\partial t}\left(  z,t\right)  -\delta_{N}^{i}%
\theta_{Ns}\dfrac{\partial^{2}C_{Ns}}{\partial z^{2}}\left(  z,t\right)   &
= & -\gamma_{is}\dfrac{\partial C_{if}}{\partial r}\left(  1,z,t\right)  \\
&  & \quad+ \delta_{i}\mathbf{r}_{i}\left(  C_{1s},\ldots,C_{Ns}\right)  \left(
z,t\right) 
\end{array}
\right.  \label{sysp}%
\end{equation}
for $i\in\left\{  1,\ldots,N\right\}  $, $\delta_{i}\in\left\{  -1,1\right\}
$ and $\delta_{N}^{i}=1$ if $i=N$, and $\delta_{N}^{i}=0$ if not, moreover the constants $\beta_{if}$, $\gamma_{is}$ and $\theta_{Ns}$ are nonnegative. The
initial and boundary conditions are
\markboth{A coupled problem modelizing an automotive catalytic converter}{Jean-David Hoernel}
\begin{equation}
\left\{
\begin{array}
[c]{rclrrcl}%
C_{if}\left(  r,0,t\right)   & = & C_{if0}\left(  r\right)  , &  &
\dfrac{\partial C_{if}}{\partial r}\left(  0,z,t\right)   & = & 0,\\
C_{if}\left(  1,z,t\right)   & = & C_{is}\left(  z,t\right)  , &  &
C_{is}\left(  z,0\right) & = &  C_{is0}\left(  z\right),\\
\theta_{Ns}\dfrac{\partial C_{Ns}}{\partial z}\left(  1,t\right)  & = & 0 , &  & \theta
_{Ns}\dfrac{\partial C_{Ns}}{\partial z}\left(  0,t\right)   & = & 0.
\end{array}
\right.  \label{cond}%
\end{equation}
The functions $\mathbf{r}_{i}$, $i\in\left\{  1,\ldots,N\right\}  $, are
supposed to verify the following assumptions
\begin{itemize}
\item[(H1)] They are Lipschitz continuous with the same constant $k$.
\item[(H2)] For all $\left(  x_{1},\ldots,x_{N}\right)  $ in $\mathbb{R}^{N}$ we have 
$\mathbf{r}_{i}\left(  x_{1},\ldots,x_{N}\right)  \geq0.$
\item[(H3)] If one of the $x_{j}$, $j\in \{1,\ldots, N\}$, is equal to $0$ then 
$\mathbf{r}_{i}\left(  x_{1},\ldots,0,\ldots,x_{N}\right)  =0.$
\end{itemize}
\begin{remark}
The physical meaning of these assumptions are
\begin{enumerate}
\item The function $\mathbf{r}_{i}$ is the speed of creation or consumption of
the $i-th$ chemical species through all the reactions in which the species appears.
\item {\rm (H2)} means that the $i$-th chemical species is either created
($\delta_{i}=1$) or consumed ($\delta_{i}=-1$).
\item {\rm (H3)} say that if there is no more of one chemical species necessary for
the reaction with the $i$-th species, then the corresponding reaction stops.
\end{enumerate}
\end{remark}
Integrating the $i$-th equation of (\ref{sysp})$_{1}$ with respect to $r$
after multiplying by $r\left(  1-r^{2}\right)  $ allows us to rewrite the system (\ref{sysp}) as
\begin{equation}
\left\{
\begin{array}
[c]{rll}%
\dfrac{\partial C_{if}}{\partial z} & = & \dfrac{\beta_{if}}{r\left(
1-r^{2}\right)  }\dfrac{\partial}{\partial r}\left(  r\dfrac{\partial C_{if}%
}{\partial r}\right)  ,\\
\dfrac{\partial C_{is}}{\partial t}-\delta_{N}^{i}\theta_{Ns}\dfrac
{\partial^{2}C_{Ns}}{\partial z^{2}} & = & -\dfrac{\gamma_{is}}{\beta_{if}}%
{\displaystyle\int\nolimits_{0}^{1}}
\dfrac{\partial C_{if}}{\partial z}r\left(  1-r^{2}\right)  dr\\
&  & +\delta_{i}\mathbf{r}_{i}\left(  C_{1s},\ldots,C_{Ns}\right)  .
\end{array}
\right.  \label{SRp}%
\end{equation}
It is this problem (\ref{SRp}) with the initials and boundary conditions (\ref{cond}) that we are studying from now.

\section{Existence of the solution}

To prove the existence of a solution of (\ref{SRp}) the main idea is to perfom 
what we called a ``parabolic regularization" of (\ref{SRp})$_{\text{2}}$ in order to avoid working 
with both pdes and odes and first show the existence of a solution of the new regularized 
problem. Then we deduce the existence of a solution of the initial problem using a limiting process.

\subsection{Regularization of the problem}

We regularized (\ref{SRp}) by adding $-\theta_{is}\frac{\partial
^{2}C_{is}}{\partial z^{2}}$ with $\theta_{is}>0$ and $i\in \{1,\ldots,N-1\}$, in the
equations (\ref{SRp})$_{\text{2}}$ on the boundary to obtain
\begin{equation}
\left\{
\begin{array}
[c]{rll}%
\dfrac{\partial C_{if}}{\partial z} & = & \dfrac{\beta_{if}}{r\left(
1-r^{2}\right)  }\dfrac{\partial}{\partial r}\left(  r\dfrac{\partial C_{if}%
}{\partial r}\right)  ,\\
\dfrac{\partial C_{is}}{\partial t}-\theta_{is}\dfrac{\partial^{2}C_{is}%
}{\partial z^{2}} & = & -\dfrac{\gamma_{is}}{\beta_{if}}%
{\displaystyle\int\nolimits_{0}^{1}}
\dfrac{\partial C_{if}}{\partial z}r\left(  1-r^{2}\right)  dr+\delta
_{i}\mathbf{r}_{i}\left(  C_{1s},\ldots,C_{Ns}\right)
\end{array}
\right.  \label{SPp_reg}%
\end{equation}
with the initials and boundary conditions (\ref{cond}) and
\[
\theta_{is}\dfrac{\partial C_{is}}{\partial z}\left(  1,t\right)
=0=\theta_{is}\dfrac{\partial C_{is}}{\partial z}\left(  0,t\right)
,\qquad\forall i\in\{1,\ldots,N-1\}.
\]
Then we prove the existence of a solution of the problem (\ref{SPp_reg}), using the
mappings
\[%
\begin{array}
[c]{ccccccc}
& \Psi &  &  &  & \Phi & \\
C_{is} & \rightarrow & C_{if} & \text{ and } & C_{if} & \rightarrow & C_{is}%
\end{array}
\]

\begin{lemma}
\label{estimthe}There exists a nonnegative constant $c\left(  T\right)  $ such
that
\begin{equation}%
{\displaystyle\sum\limits_{i=1}^{N}}
{\displaystyle\int\nolimits_{0}^{1}}
\left(  C_{is}\right)  ^{2}\left(  z,T\right)  dz+%
{\displaystyle\sum\limits_{i=1}^{N}}
\theta_{is}%
{\displaystyle\int\nolimits_{0}^{T}}
{\displaystyle\int\nolimits_{0}^{1}}
\left(  \dfrac{\partial C_{is}}{\partial z}\right)  ^{2}dzdt\leq c\left(
T\right)  . \label{ineg_theta_is}%
\end{equation}
\end{lemma}
\vspace{-3mm}
\begin{proof}
From {\rm (\ref{SPp_reg})}, we deduce after summation on $i$ and using {\rm (H1)} and {\rm (H3)} that
\[%
\begin{array}
[c]{l}%
\dfrac{1}{2}%
{\displaystyle\sum\limits_{i=1}^{N}}
\dfrac{\beta_{if}}{\gamma_{is}}%
{\displaystyle\int\nolimits_{0}^{l}}
\left(  C_{is}\right)  ^{2}\left(  z,T\right)  dz+%
{\displaystyle\sum\limits_{i=1}^{N}}
\theta_{is}\dfrac{\beta_{if}}{\gamma_{is}}%
{\displaystyle\int\nolimits_{0}^{T}}
{\displaystyle\int\nolimits_{0}^{l}}
\left(  \dfrac{\partial C_{is}}{\partial z}\right)  ^{2}dzdt\\
\leq
a(T)+\underset{i}{\sup}\left(  \dfrac{\beta_{if}}{\gamma_{is}}\right)
kN%
{\displaystyle\sum\limits_{i=1}^{N}}
{\displaystyle\int\nolimits_{0}^{T}}
{\displaystyle\int\nolimits_{0}^{l}}
\left(  C_{is}\right)  ^{2}dzdt.
\end{array}
\]
\vspace{-3mm}

\noindent with $a\left(  T\right)  =\dfrac{T}{2}%
{\displaystyle\sum\limits_{i=1}^{N}}
{\displaystyle\int\nolimits_{0}^{1}}
\left(  C_{i0}\right)  ^{2}\left(  r\right)  r\left(  1-r^{2}\right)
dr+\dfrac{1}{2}%
{\displaystyle\sum\limits_{i=1}^{N}}
\dfrac{\beta_{if}}{\gamma_{is}}%
{\displaystyle\int\nolimits_{0}^{1}}
\left(  C_{is0}\right)  ^{2}\left(  z\right)  dz,$ which leads to%
\[%
{\displaystyle\sum\limits_{i=1}^{N}}
{\displaystyle\int\nolimits_{0}^{l}}
\left(  C_{is}\right)  ^{2}\left(  z,T\right)  dz\leq b(T)+d%
{\displaystyle\sum\limits_{i=1}^{N}}
{\displaystyle\int\nolimits_{0}^{T}}
{\displaystyle\int\nolimits_{0}^{l}}
\left(  C_{is}\right)  ^{2}dzdt
\]
with $b(T)=2a(T)/\inf\nolimits_{i}\left(  \frac{\beta_{if}}%
{\gamma_{is}}\right)  $ and $d=2\sup\nolimits_{i}\left(
\frac{\beta_{if}}{\gamma_{is}}\right)  kN/\inf\nolimits_{i}\left(  \frac
{\beta_{if}}{\gamma_{is}}\right)  .$ Using Gronwall's lemma leads to the result.
\end{proof}

\subsection{Existence in the cylinder}

Assume that $C_{is}$, $i\in\{1,\ldots,N\}$ are known on the boundary (mapping $\Psi$) and let
\[
\left\{
\begin{array}
[c]{rrllll}%
u_{f} & = & ^{t}\left(  C_{1f},\ldots,C_{Nf}\right)  , & u_{f0} & = &
^{t}\left(  C_{1f0},\ldots,C_{Nf0}\right)  ,\\
u_{s} & = & ^{t}\left(  C_{1s},\ldots,C_{Ns}\right)  , & \beta_{f} & = &
diag\left(  \beta_{1f},\cdots,\beta_{Nf}\right),
\end{array}
\right.
\]
to obtain the following problem in the cylinder
\begin{equation}
\left\{
\begin{array}
[c]{rll}%
\dfrac{\partial u_{f}}{\partial z}-\beta_{f}\dfrac{1}{r\left(  1-r^{2}\right)
}\dfrac{\partial}{\partial r}\left(  r\dfrac{\partial u_{f}}{\partial
r}\right)  & = & 0,\\
u_{f}\left(  r,0,t\right)  & = & u_{f0}\left(  r\right)  ,\\
u_{f}\left(  1,z,t\right)  & = & u_{s}\left(  z,t\right)  ,\\
\dfrac{\partial u_{f}}{\partial r}\left(  0,z,t\right)  & = & 0.
\end{array}
\right.  \label{Suf}%
\end{equation}
Performing the change of fonction $w_{f}\left(  r,z,t\right)  =u_{f}\left(  r,z,t\right)  -u_{s}\left(
z,t\right)$ we have the problem
\begin{equation}
\left\{
\begin{array}
[c]{rll}%
\dfrac{\partial w_{f}}{\partial z}-\beta_{f}\dfrac{1}{r\left(  1-r^{2}\right)
}\dfrac{\partial}{\partial r}\left(  r\dfrac{\partial w_{f}}{\partial
r}\right)  & = & -\dfrac{\partial u_{s}}{\partial z},\\
w_{f}\left(  r,0,t\right)  & = & w_{0}\left(  r,t\right)  ,\\
w_{f}\left(  1,z,t\right)  & = & 0,\\
\dfrac{\partial w_{f}}{\partial r}\left(  0,z,t\right)  & = & 0,
\end{array}
\right.  \label{SH}%
\end{equation}
with $w_{0}\left(  r,t\right)  =u_{f0}\left(  r\right)  -u_{s}\left(
0,t\right)  =u_{f0}\left(  r\right)  -u_{f0}\left(  1\right)  $.

\begin{definition}
\begin{enumerate}
\item We set
$$
\begin{tabular}
[c]{rll}%
$L_{r(1-r^{2})}^{2}\left(  0,1\right)  $ & $=$ & $\left\{  u\mid r\mapsto
u\left(  r\right)  \sqrt{r\left(  1-r^{2}\right)  }\in L^{2}\left(
0,1\right)  \right\}  ,$\\
$L_{r}^{2}\left(  0,1\right)  $ & $=$ & $\left\{  u\mid r\mapsto u\left(
r\right)  \sqrt{r}\in L^{2}\left(  0,1\right)  \right\}  .$%
\end{tabular}
$$
$L_{r(1-r^{2})}^{2}\left(  0,1\right)  $ (resp. $L_{r}^{2}\left(  0,1\right)
$) is a Hilbert space for the norm associated to the scalar product
\vspace{-3mm}

$$\left\langle u,v\right\rangle _{2,r(1-r^{2})}=\int\nolimits_{0}^{1}uvr\left(
1-r^{2}\right)  dr\text{ (resp. }\left\langle u,v\right\rangle _{2,r}%
=\int\nolimits_{0}^{1}uvrdr\text{)}.$$
\vspace{-5mm}

\item Let
$$W_{r}=\left\{  u\in\left(  L_{r(1-r^{2})}^{2}\left(  0,1\right)  \right)
^{N}\mid\dfrac{\partial u}{\partial r}\in\left(  L_{r}^{2}\left(  0,1\right)
\right)  ^{N}\right\} ,$$
$$W_{r0}=\left\{  u\in W_{r}\mid u\left(  1\right)  =0\right\}  ,\text{ }%
W_{r}\left(  T\right)  =\left\{  u\in L^{2}\left(  \left]  0,1\right[
\times\left]  0,T\right[  ;W_{r}\right)  \right\}$$
and $W_{r0}^{\prime}$ be the dual space of $W_{r0}$.

\item The spaces $W_{r}$, $W_{r0}$ and $W_{r}\left(  T\right)  $ are equipped
with the norms
$$
\left\Vert u\right\Vert _{W_{r}}^{2}:=%
{\displaystyle\int\nolimits_{0}^{1}}
\left\Vert u\right\Vert ^{2}r\left(  1-r^{2}\right)  dr+%
{\displaystyle\int\nolimits_{0}^{1}}
\left\Vert \dfrac{\partial u}{\partial r}\right\Vert ^{2}rdr,\text{
}\left\Vert u\right\Vert _{W_{r0}}^{2}:=%
{\displaystyle\int\nolimits_{0}^{1}}
\left\Vert \dfrac{\partial u}{\partial r}\right\Vert ^{2}rdr,
$$
$$
\left\Vert u\right\Vert _{W_{r}(T)}^{2}:=%
{\displaystyle\int\nolimits_{0}^{T}}
{\displaystyle\int\nolimits_{0}^{1}}
{\displaystyle\int\nolimits_{0}^{1}}
\left\Vert u\right\Vert ^{2}r\left(  1-r^{2}\right)  drdzdt+%
{\displaystyle\int\nolimits_{0}^{T}}
{\displaystyle\int\nolimits_{0}^{1}}
{\displaystyle\int\nolimits_{0}^{1}}
\left\Vert \dfrac{\partial u}{\partial r}\right\Vert ^{2}rdrdzdt,
$$
where $\left\Vert .\right\Vert $ is the Euclidian norm of $\mathbb{R}^{N}$.
\end{enumerate}
\end{definition}
We establish some properties of these spaces
\begin{lemma}
\begin{enumerate}
\item We have the continuous embeddings
$$L^{2}\left(  0,1\right)  \underset{\rightarrow}{\subset}L_{r}^{2}\left(
0,1\right)  \underset{\rightarrow}{\subset}L_{r(1-r^{2})}^{2}\left(
0,1\right)  .$$ 
\vspace{-8mm}
\item For all $u$ in $W_{r0}$, we have
\begin{equation}%
{\displaystyle\int\nolimits_{0}^{1}}
\left\Vert u\right\Vert ^{2}r\left(  1-r^{2}\right)  dr\leq\dfrac{3}{16}%
{\displaystyle\int\nolimits_{0}^{1}}
\left\Vert \dfrac{\partial u}{\partial r}\right\Vert ^{2}rdr. \label{pf}%
\end{equation}
Then the injection of $W_{r0}$ in $\left(  L_{r(1-r^{2})}%
^{2}\left(  0,1\right)  \right)  ^{N}$ is continuous and compact and the embedding is dense.
\item $W_{r0}$ is a Hilbert space for the norm above.
\end{enumerate}
\end{lemma}
\begin{proof}
1. As $\forall r\in\left[  0,1\right] ,$ $0\leq
r\left(  1-r^{2}\right)  \leq r\leq1,$ we deduce the result from Kufner {\rm [9] .}

\noindent2. For all $a$ such that $0<a<1$ we have as $u(1)=0$ that
\[
\left\Vert u(a)\right\Vert =\left\Vert -%
{\displaystyle\int\nolimits_{a}^{1}}
\dfrac{\partial u}{\partial r}dr\right\Vert \leq%
{\displaystyle\int\nolimits_{a}^{1}}
\left\Vert \dfrac{\partial u}{\partial r}\right\Vert \sqrt{r}\dfrac{1}%
{\sqrt{r}}dr.
\]
Using the Cauchy-Schwarz inequalitie, Integrating and multiplying by $a(1-a^{2})$ leads to the result.
We show the compactness of the injection of $W_{r0}$ in $\left(
L_{r(1-r^{2})}^{2}\left(  0,1\right)  \right)  ^{N}$ and the density of the embedding as in Dautray-Lions {\rm [4].} For more details, see Hoernel {\rm [6].}

\noindent3. Immediate consequence of 2.
\end{proof}

\begin{definition}
\label{defcyl}Let $u_{f0}\in\left(  L_{r(1-r^{2})}^{2}\left(  0,1\right)  \right)  ^{N}$ and
$u_{s}\in L^{2}\left(  0,T;\left(  H^{1}\left(  0,1\right)  \right)
^{N}\right).$ A function $w_{f}$ is called a \textbf{weak solution} of {\rm (\ref{SH})} if and
only if $w_{f}\in W_{r}\left(  T\right),$ 
$\frac{\partial w_{f}}{\partial z}\in L^{2}\left(  \left]  0,1\right[  \times\left]  0,T\right[
;W_{r0}^{\prime}\right)  $, $w_{f}\left(  r,0,t\right)  =w_{0}\left(
r,t\right)  $, and if for every test-function $\varphi\in L^{2}\left(  \left]
0,1\right[  \times\left]  0,T\right[  ;W_{r0}\right)  $ we have
\begin{equation}%
\begin{array}
[c]{l}%
{\displaystyle\int\nolimits_{0}^{T}}
{\displaystyle\int\nolimits_{0}^{1}}
{\displaystyle\int\nolimits_{0}^{1}}
\left(  \dfrac{\partial w_{f}}{\partial z}\cdot\varphi\right)  r\left(
1-r^{2}\right)  drdzdt+%
{\displaystyle\int\nolimits_{0}^{T}}
{\displaystyle\int\nolimits_{0}^{1}}
{\displaystyle\int\nolimits_{0}^{1}}
\left(  \beta_{f}\dfrac{\partial w_{f}}{\partial r}\cdot\dfrac{\partial
\varphi}{\partial r}\right)  rdrdzdt\\
=-%
{\displaystyle\int\nolimits_{0}^{T}}
{\displaystyle\int\nolimits_{0}^{1}}
\left(  \dfrac{\partial u_{s}}{\partial z}\cdot%
{\displaystyle\int\nolimits_{0}^{1}}
\varphi r\left(  1-r^{2}\right)  dr\right)  dzdt.
\end{array}
\label{formvarf}%
\end{equation}
\end{definition}
\vspace{-4mm}
To prove the existence of a weak solution of (\ref{SH}) we define the second
order elliptic operator $L:W_{r0}\rightarrow W_{r0}^{\prime}$ by
\vspace{-3mm}

\[
\forall w,v\in W_{r0}:\left\langle Lw,v\right\rangle _{\left\langle
W_{r0}^{\prime},W_{r0}\right\rangle }=%
{\displaystyle\int_{0}^{1}}
\left(  \beta_{f}\dfrac{\partial w}{\partial r}\cdot\dfrac{\partial
v}{\partial r}\right)  rdr
\]
\vspace{-2mm}

\noindent and the operator $T$ from $\left(  L_{r(1-r^{2})}^{2}\left(
0,1\right)  \right)  ^{N}$ in itself giving for each function $g$ taken in $\left(
L_{r(1-r^{2})}^{2}\left(  0,1\right)  \right)  ^{N}$ the unique element
$T\left(  g\right)  $ in $W_{r0}$ that is the weak solution of $L\left(
T\left(  g\right)  \right)  =gr\left(  1-r^{2}\right)  $ by
\begin{equation}%
{\displaystyle\int\nolimits_{0}^{1}}
\left(  \beta_{f}\dfrac{\partial T\left(  g\right)  }{\partial r}\cdot
\dfrac{\partial\varphi}{\partial r}\right)  rdr=%
{\displaystyle\int\nolimits_{0}^{1}}
\left(  g\cdot\varphi\right)  r\left(  1-r^{2}\right)  dr\text{,\qquad}%
\forall\varphi\in W_{r0}. \label{LT}%
\end{equation}

\begin{lemma}
\begin{enumerate}
\item The operator $T$ is

- Self-adjoint: $\forall f,g\in\left(  L_{r(1-r^{2})}^{2}\left(
0,1\right)  \right)  ^{N}:
\left\langle T\left(  f\right)  ,g\right\rangle _{r(1-r^{2})}=\left\langle
f,T\left(  g\right)  \right\rangle _{r(1-r^{2})}$

- Positive: $\forall f\in\left(  L_{r(1-r^{2})}^{2}\left(  0,1\right)
\right)  ^{N}:\left\langle T\left(  f\right)  ,f\right\rangle _{r(1-r^{2}%
)}\geq0$

- Non-degenerated: $\forall f\in\left(  L_{r(1-r^{2})}^{2}\left(
0,1\right)  \right)  ^{N}:\left\langle T\left(  f\right)  ,f\right\rangle
_{r(1-r^{2})}=0\Rightarrow f=0.$
\item Moreover there exists a Hilbert basis $\left\{  \omega^{j}\right\}
_{j}$ of both $\left(  L_{r(1-r^{2})}^{2}\left(  0,1\right)  \right)  ^{N}$
and $W_{r0}$ made of eigenfunctions of $T$ and $L$.
\end{enumerate}
\end{lemma}

\begin{proof}
1. Immediate.

\noindent2. As the injection of $W_{r0}$ in $\left(  L_{r(1-r^{2})}^{2}\left(
0,1\right)  \right)  ^{N}$ is compact, $T$ is a compact operator from $\left(
L_{r(1-r^{2})}^{2}\left(  0,1\right)  \right)  ^{N}$ in itself. Taking into
acount the properties of $T$, the Theorem {\rm VI.11} from Brezis {\rm [1, page 97]}
implies the existence of a Hilbert basis $\left\{  \omega^{j}\right\}  _{j}$
of $\left(  L_{r(1-r^{2})}^{2}\left(  0,1\right)  \right)  ^{N}$ made of
eigenfunctions of $T$\ and of $L$.\ The Remark {\rm 29} from Brezis {\rm [1, page 193]}
implies that this basis is also a Hilbert basis of $W_{r0}$.
\end{proof}

\noindent We show the following existence result
\begin{proposition}
\label{solfai}Let $u_{0}$ and $u_{s}$ as in Definition {\rm \ref{defcyl}}. Then,
there exists at least one \textbf{weak solution }$w_{f}$ of {\rm (\ref{SH})}.
\end{proposition}

\begin{proof}
Let $\left\{  \omega^{j}\right\}  _{j}$ be a basis of $W_{r0}$ made of
eigenfunctions of $L$ with eigenvalues $\lambda_{j}$
\[
\left\langle L\left(  \omega^{j}\right)  ,\varphi\right\rangle _{\left\langle
W_{r0}^{\prime},W_{r0}\right\rangle }=\lambda_{j}%
{\displaystyle\int\nolimits_{0}^{1}}
\left(  \omega^{j}\cdot\varphi\right)  r\left(  1-r^{2}\right)  dr,\qquad
\forall\varphi\in W_{r0}.
\]
For all $\varphi$ in $W_{r0}$, there exists a sequence $\left(  \gamma
_{j}\right)  _{j}$ such that $\varphi=\sum_{j=1}^{+\infty}\gamma_{j}\omega
^{j}$. We write $P_{m}\left(  \varphi\right)  =\sum_{j=1}^{m}\gamma_{j}%
\omega^{j}$ the Galerkin approximation of order $m$ of $\varphi$. Properties of
$\left\{  \omega^{j}\right\}  _{j}$ implies that the operator $P_{m}$ is a
continuous orthogonal projector from $W_{r0}$ (resp. \newline $\left(  L_{r(1-r^{2}%
)}^{2}\left(  0,1\right)  \right)  ^{N}$) to the subspace of $W_{r0}$ (resp.
of $\left(  L_{r(1-r^{2})}^{2}\left(  0,1\right)  \right)  ^{N}$) generated by
$\omega^{1},\ldots,\omega^{m}$. Moreover
\[
\left\Vert P_{m}\left(  \varphi\right)  \right\Vert _{W_{r0}}\leq\left\Vert
\varphi\right\Vert _{W_{r0}}\text{ ; }\left\Vert P_{m}\left(  \varphi\right)
\right\Vert _{2,r(1-r^{2})}\leq\left\Vert \varphi\right\Vert _{2,r(1-r^{2}%
)},\text{\qquad}\forall\varphi\in W_{r0}.
\]
We establish the following estimates
\begin{equation}%
\begin{array}
[c]{rrl}%
{\displaystyle\int\nolimits_{0}^{T}}
{\displaystyle\int\nolimits_{0}^{1}}
\left\Vert P_{m}\left(  w_{f}\right)  \right\Vert ^{2}\left(  r,s,t\right)
r(1-r^{2})drdt & \leq & C(T),\\%
{\displaystyle\int\nolimits_{0}^{T}}
{\displaystyle\int\nolimits_{0}^{s}}
{\displaystyle\int\nolimits_{0}^{1}}
\left\Vert \dfrac{\partial P_{m}\left(  w_{f}\right)  }{\partial r}\right\Vert
^{2}\left(  r,z,t\right)  rdrdzdt & \leq & C(T),\\
\left\Vert \dfrac{\partial P_{m}\left(  w_{f}\right)  }{\partial z}\right\Vert
_{L^{2}\left(  \left]  0,1\right[  \times\left]  0,T\right[  ;W_{r0}^{\prime
}\right)  } & \leq & C(T),
\end{array}
\label{estm}%
\end{equation}
for almost every $s$ in $\left]  0,1\right[  $, where $C(T)$ is a nonnegative 
constant that only depends of the problem {\rm (\ref{SH})} datas and $T$,
but independant of $m$. Then, there
exists a subsequence of $\left(  P_{m}\left(  w_{f}\right)  \right)  _{m}$
still noted $\left(  P_{m}\left(  w_{f}\right)  \right)  _{m}$ such that
\[%
\begin{array}
[c]{rrl}%
P_{m}\left(  w_{f}\right)   & \underset{m\rightarrow+\infty}{\rightharpoonup}
& w_{f}\text{, weakly in }L^{2}\left(  \left]  0,1\right[  \times\left]
0,T\right[  ;\left(  L_{r(1-r^{2})}^{2}(0,1)\right)  ^{N}\right)  ,\\
\dfrac{\partial P_{m}\left(  w_{f}\right)  }{\partial r} & \underset
{m\rightarrow+\infty}{\rightharpoonup} & \frac{\partial w_{f}}{\partial r}\text{, weakly in }%
L^{2}\left(  \left]  0,1\right[  \times\left]  0,T\right[  ;\left(  L_{r}%
^{2}(0,1)\right)  ^{N}\right)  ,\\
\dfrac{\partial P_{m}\left(  w_{f}\right)  }{\partial z} & \underset
{m\rightarrow+\infty}{\rightharpoonup} & \frac{\partial w_{f}}{\partial z}\text{, weakly in }%
L^{2}\left(  \left]  0,1\right[  \times\left]  0,T\right[  ;W_{r0}^{\prime
}\right)  .
\end{array}
\]
Consequently $\left(  P_{m}\left(  w_{f}\right)  \right)
_{m}$ weakly converges to $w_{f}$ in $L^{2}\left(  \left]  0,1\right[
\times\left]  0,T\right[  ;W_{r0}\right)  $ and  $w_{f}$ verifies the variationnal
formulation {\rm (\ref{formvarf})}.
\end{proof}
\begin{remark}
If $u_{f0}$ belongs to $W_{r}$ then
$
\frac{\partial w_{f}}{\partial z}\in L^{2}\left(  \left]  0,1\right[
\times\left]  0,T\right[  ;\left(  L_{r\left(  1-r^{2}\right)  }^{2}\left(
0,1\right)  \right)  ^{N}\right)
$
\end{remark}
\vspace{-4mm}

\begin{remark}
From
$L^{2}\left(  \left]  0,1\right[  \times\left]  0,T\right[  ;W_{r0}\right)
\subset L^{2}\left(  \left]  0,1\right[  \times\left]  0,T\right[  ;\left(
L_{r(1-r^{2})}^{2}\left(  0,1\right)  \right)  ^{N}\right)  ,$\newline
$L^{2}\left(  \left]  0,1\right[  \times\left]  0,T\right[  ;\left(
L_{r(1-r^{2})}^{2}\left(  0,1\right)  \right)  ^{N}\right)  \subset
L^{2}\left(  \left]  0,1\right[  \times\left]  0,T\right[  ;W_{r0}^{\prime
}\right)$
and the fact that we have \newline $w_{f}\in L^{2}\left(  \left]  0,1\right[  \times\left]
0,T\right[  ;W_{r0}\right)  $ and $\frac{\partial w_{f}}{\partial z}\in
L^{2}\left(  \left]  0,1\right[  \times\left]  0,T\right[  ;W_{r0}^{\prime
}\right)  $, using Proposition {\rm 23.23} from Zeidler {\rm [13, page 422]} we 
obtain that  $w_{f}\in C\left(  \left]  0,1\right[  \times\left]
0,T\right[  ;\left(  L_{r(1-r^{2})}^{2}\left(  0,1\right)  \right)
^{N}\right)  $ which allows us to interpret the initial condition.
\end{remark}
\vspace{-3mm}

\begin{corollary}
The problem {\rm (\ref{Suf})} admits a weak solution $u_{f}\in W_{r}\left(  T\right)
$ such that
\[%
{\displaystyle\int\nolimits_{0}^{T}}
{\displaystyle\int\nolimits_{0}^{1}}
{\displaystyle\int\nolimits_{0}^{1}}
\left(  \dfrac{\partial u_{f}}{\partial z}\cdot\varphi\right)  r\left(
1-r^{2}\right)  drdzdt+%
{\displaystyle\int\nolimits_{0}^{T}}
{\displaystyle\int\nolimits_{0}^{1}}
{\displaystyle\int\nolimits_{0}^{1}}
\left(  \beta_{f}\dfrac{\partial u_{f}}{\partial r}\cdot\dfrac{\partial
\varphi}{\partial r}\right)  rdrdzdt=0,
\]
for all $T>0$ and all $\varphi\in L^{2}\left(  \left]  0,1\right[
\times\left]  0,T\right[  ;W_{r0}\right)  $ and $u_{f}$ verify $u_{f}\left(
r,0,t\right)  =u_{f0}\left(  r\right)  $.
\end{corollary}

\begin{proof}
Let us notice that $u_{f}=w_{f}+u_{s}$ is a weak solution of {\rm (\ref{Suf})}
belonging to the right space and verifying the initial condition at $z=0$.
\end{proof}

\begin{remark}
\label{dualite}Since $\frac{\partial u_{f}}{\partial z}\in L^{2}\left(
\left]  0,1\right[  \times\left]  0,T\right[  ;W_{r0}^{\prime}\right)  $ and
the function $r\mapsto r\left(  1-r^{2}\right)  $ is equal to zero at $r=1$,
we can affirm that for all function $g:$ $\left(  z,t\right)  \mapsto g\left(
z,t\right)  $ in $\left(  L^{2}\left(  \left]  0,1\right[  \times\left]
0,T\right[  \right)  \right)  ^{N}$, the function $\left(  r,z,t\right)
\mapsto g\left(  z,t\right)  r\left(  1-r^{2}\right)  $ is in the space
$L^{2}\left(  \left]  0,1\right[  \times\left]  0,T\right[  ;W_{r0}\right)  $
and the following duality bracket is well defined
\[
\left\langle \dfrac{\partial u_{f}}{\partial z},gr\left(  1-r^{2}\right)
\right\rangle _{\left\langle L^{2}\left(  \left]  0,1\right[  \times\left]
0,T\right[  ;W_{r0}^{\prime}\right)  ,L^{2}\left(  \left]  0,1\right[
\times\left]  0,T\right[  ;W_{r0}\right)  \right\rangle }%
\]
\end{remark}
\vspace{-2mm}

\subsection{Existence on the boundary}

We use the mapping $\Phi$. Suppose given the $\frac{\partial C_{if}}{\partial z}$ for $i\in\{1,\ldots,N\}$ and let
\[
\left\{
\begin{array}
[c]{rlrl}%
u_{s0} & =^{t}\left(  C_{1s0},\ldots,C_{Ns0}\right)  , & \delta & =diag\left(
\delta_{1},\cdots,\delta_{N}\right)  ,\\
\mathbf{r} & =^{t}\left(  \mathbf{r}_{1,}\ldots,\mathbf{r}_{N}\right)  , &
\theta_{s} & =diag\left(  \theta_{1s},\cdots,\theta_{Ns}\right)  ,\\
\Gamma & =diag\left(  \dfrac{\gamma_{1s}}{\beta_{1f}},\cdots,\dfrac
{\gamma_{Ns}}{\beta_{Nf}}\right)  & =: & diag\left(  \Gamma_{1},\ldots
,\Gamma_{N}\right)  ,
\end{array}
\right.
\]
in order to obtain the system
\begin{equation}
\left\{
\begin{array}
[c]{rll}%
\dfrac{\partial u_{s}}{\partial t}-\theta_{s}\dfrac{\partial^{2}u_{s}%
}{\partial z^{2}} & = & \delta\mathbf{r}(u_{s})-\Gamma%
{\displaystyle\int\nolimits_{0}^{1}}
\dfrac{\partial u_{f}}{\partial z}r\left(  1-r^{2}\right)  dr,\\
u_{s}\left(  z,0\right)  & = & u_{s0}\left(  z\right)  ,\\
\theta_{s}\dfrac{\partial u_{s}}{\partial z}\left(  0,t\right)  & = &
0=\theta_{s}\dfrac{\partial u_{s}}{\partial z}\left(  1,t\right)  .
\end{array}
\right.  \label{sysbord}%
\end{equation}
\vspace{-3mm}

\noindent Let
\vspace{-2mm}
\[
W_{z}=\left\{  u\in\left(  L^{2}\left(  0,1\right)  \right)  ^{N}\mid
\dfrac{\partial u}{\partial z}\in\left(  L^{2}\left(  0,1\right)  \right)
^{N}\right\}  ,\text{ }W_{z}\left(  T\right)  =\left\{  u\in L^{2}\left(
0,T;W_{z}\right)  \right\} 
\]
\vspace{-2mm}
and $W_{z}^{\prime}$ be the dual space of $W_{z}$. These spaces are equipped
with the following norms
\begin{align*}
\left\Vert u\right\Vert _{W_{z}}^{2}  & :=%
{\displaystyle\int\nolimits_{0}^{1}}
\left\Vert u\right\Vert ^{2}dz+%
{\displaystyle\int\nolimits_{0}^{1}}
\left\Vert \dfrac{\partial u}{\partial z}\right\Vert ^{2}dz,\text{ }\\
\left\Vert u\right\Vert _{W_{z}(T)}^{2}  & :=%
{\displaystyle\int\nolimits_{0}^{T}}
{\displaystyle\int\nolimits_{0}^{1}}
\left\Vert u\right\Vert ^{2}dzdt+%
{\displaystyle\int\nolimits_{0}^{T}}
{\displaystyle\int\nolimits_{0}^{1}}
\left\Vert \dfrac{\partial u}{\partial z}\right\Vert ^{2}dzdt.
\end{align*}
\vspace{-4mm}

\begin{definition}
\label{solfaip}Assume $
\frac{\partial u_{f}}{\partial z}\in L^{2}\left(  \left]  0,T\right[
\times\left]  0,1\right[  ;W_{r0}^{\prime}\right)$ and $u_{s0}\in\left(
L^{2}\left(  0,1\right)  \right)  ^{N}.$ A function $u_{s}$ is called a \textbf{weak solution} of {\rm (\ref{sysbord})} if
and only if we have $u_{s}\in W_{z}\left(  T\right)  $, $u_{s}\left(  z,0\right)
=u_{s0}\left(  z\right)  $, and if, for all $\psi$ in $L^{2}\left(
0,T;W_{z}\right)  $, we have
\begin{equation}%
\begin{array}
[c]{l}%
{\displaystyle\int\nolimits_{0}^{\tau}}
{\displaystyle\int\nolimits_{0}^{1}}
\left(  \dfrac{\partial u_{s}}{\partial t}\cdot\psi\right)  dzdt+%
{\displaystyle\int\nolimits_{0}^{\tau}}
{\displaystyle\int\nolimits_{0}^{1}}
\left(  \theta_{s}\dfrac{\partial u_{s}}{\partial z}\cdot\dfrac{\partial\psi
}{\partial z}\right)  dzdt\\
=%
{\displaystyle\int\nolimits_{0}^{\tau}}
{\displaystyle\int\nolimits_{0}^{1}}
\left(  \delta\mathbf{r}(u_{s})\cdot\psi\right)  dzdt-%
{\displaystyle\int\nolimits_{0}^{\tau}}
{\displaystyle\int\nolimits_{0}^{1}}
\left\langle \Gamma\dfrac{\partial u_{f}}{\partial z},\psi r\left(
1-r^{2}\right)  \right\rangle _{\left\langle W_{r0}^{\prime},W_{r0}%
\right\rangle }dzdt,
\end{array}
\label{formz}%
\end{equation}
\vspace{-2mm}

\noindent for $\tau\leq T$, cf. Remark {\rm \ref{dualite}}.
\end{definition}
\vspace{-4mm}
\begin{proposition}
\label{existpar}Let $u_{f}$ and $u_{s0}$ be as in the Definition
{\rm \ref{solfaip}}.\ Then, there exists at least one \textbf{weak solution }$u_{s}$
of {\rm (\ref{sysbord})}.
\end{proposition}
\begin{proof}
We use some auxilliary linearized equation and some fixed point as in Chipot {\rm [3]} or Henry {\rm [5]}. For more details, see Hoernel {\rm [6]}.
\end{proof}
%

\vspace{2mm}
\subsection{Existence}

We begin by showing some properties of $\Phi$ and $\Psi$
\begin{proposition}
The mapping
$$\Phi:\left(
\begin{array}
[c]{ccc}%
W_{r}\left(  T\right)  & \rightarrow & W_{z}\left(  T\right) \\
u_{f} & \mapsto & u_{s}%
\end{array}
\right)$$
which for every $u_{f}\in W_{r}\left(  T\right)  $ gives the weak solution
$u_{s}\in W_{z}\left(  T\right)  $ of {\rm (\ref{sysbord})} verifies
\begin{equation}
\left\Vert \Phi\left(  u_{f}^{1}\right)  -\Phi\left(  u_{f}^{2}\right)
\right\Vert _{W_{z}\left(  T\right)  }^{2}\leq T\left(  aT+b\right)
{\displaystyle\int\nolimits_{0}^{T}}
{\displaystyle\int\nolimits_{0}^{1}}
\left\Vert \dfrac{\partial\left(  u_{f}^{1}-u_{f}^{2}\right)  }{\partial
z}\right\Vert _{W_{r0}^{\prime}}^{2}dzdt \label{phi}%
\end{equation}
where $a$ and $b$ are two nonnegative constants independant on $T$.
\end{proposition}

\begin{proof}
Let $U_{s}=u_{s}^{1}-u_{s}^{2}$ be the difference of two solutions of
{\rm (\ref{sysbord})} with same initial conditions at $t=0$ associated to $u_{f}^{1}$ and $u_{f}^{2}$ in $W_{r}\left(  T\right)  $.\ $U_{s}$ is a weak solution of
\begin{equation}
\left\{
\begin{array}
[c]{rll}%
\dfrac{\partial U_{s}}{\partial t}-\theta_{s}\dfrac{\partial^{2}U_{s}%
}{\partial z^{2}} & = & -%
{\displaystyle\int\nolimits_{0}^{1}}
\Gamma\dfrac{\partial U_{f}}{\partial z}r\left(  1-r^{2}\right)
dr+\delta\left(  \mathbf{r}(u_{s}^{1})-\mathbf{r}\left(  u_{s}^{2}\right)
\right)  ,\\
U_{s}\left(  z,0\right)  & = & 0,\\
\theta_{s}\dfrac{\partial U_{s}}{\partial z}\left(  0,t\right)  & = &
0=\theta_{s}\dfrac{\partial U_{s}}{\partial z}\left(  1,t\right)  ,
\end{array}
\right.  \label{syslipbord}%
\end{equation}
with $U_{f}=u_{f}^{1}-u_{f}^{2}$. Multiplying {\rm (\ref{syslipbord})$_\text{1}$} by
$U_{s},$ integrating on $\left]  0,1\right[  \times\left]  0,\tau\right[  $
with $0\leq\tau\leq T$ and using {\rm (H1)} leads to
\[%
\begin{array}
[c]{l}%
\dfrac{1}{2}%
{\displaystyle\int\nolimits_{0}^{1}}
\left\Vert U_{s}\right\Vert ^{2}(z,\tau)dz+%
{\displaystyle\int\nolimits_{0}^{\tau}}
{\displaystyle\int\nolimits_{0}^{1}}
\left\Vert \theta_{s}^{1/2}\dfrac{\partial U_{s}}{\partial z}\right\Vert
^{2}dzdt\\
\leq k%
{\displaystyle\int\nolimits_{0}^{\tau}}
{\displaystyle\int\nolimits_{0}^{1}}
\left\Vert U_{s}\right\Vert ^{2}dzdt+c%
{\displaystyle\int\nolimits_{0}^{\tau}}
{\displaystyle\int\nolimits_{0}^{1}}
\left\Vert \Gamma\dfrac{\partial U_{f}}{\partial z}\right\Vert _{W_{r0}%
^{\prime}}\left\Vert U_{s}\right\Vert dzdt
\end{array}
\]
with $c>0.$ Using Young's inequality with $\varepsilon>0$ and Gronwall's lemma we have
\[%
\begin{array}
[c]{c}
{\displaystyle\int\nolimits_{0}^{T}}
{\displaystyle\int\nolimits_{0}^{1}}
\left\Vert \dfrac{\partial U_{s}}{\partial z}\right\Vert ^{2}dzdt\leq
\dfrac{c\varepsilon e^{(c+2k)T/\varepsilon}}{2\inf_{i}\theta_{is}}%
{\displaystyle\int\nolimits_{0}^{T}}
{\displaystyle\int\nolimits_{0}^{1}}
\left\Vert \Gamma\dfrac{\partial U_{f}}{\partial z}\right\Vert _{W_{r0}%
^{\prime}}^{2}dzdt
\end{array}
\]
and taking $\varepsilon=(c+2k)T$ we finally obtain
\[
\left\Vert U_{s}\right\Vert _{W_{z}\left(  T\right)  }^{2}\leq\left(
c(c+2k)\left(  e-1\right)  T^{2}+\dfrac{c(c+2k)e}{2\inf_{i}\theta_{is}%
}T\right)
{\displaystyle\int\nolimits_{0}^{T}}
{\displaystyle\int\nolimits_{0}^{1}}
\left\Vert \Gamma\dfrac{\partial U_{f}}{\partial z}\right\Vert _{W_{r0}%
^{\prime}}^{2}dzdt.
\]
\end{proof}

\begin{proposition}
\label{gamma}The mapping
\[
\Psi:\left(
\begin{array}
[c]{ccc}%
W_{z}\left(  T\right)  & \rightarrow & W_{r}\left(  T\right) \\
u_{s} & \mapsto & u_{f}%
\end{array}
\right)
\]
which for every $u_{s}\in$ $W_{z}\left(  T\right)  $ gives the weak solution
$u_{f}$ of {\rm (\ref{Suf})} is such that %
\begin{equation}
\left\Vert \Psi\left(  u_{s}^{1}\right)  -\Psi\left(  u_{s}^{2}\right)
\right\Vert _{W_{r}(T)}\leq c\left\Vert u_{s}^{1}-u_{s}^{2}\right\Vert
_{W_{z}(T)} \label{fluid}%
\end{equation}
for a nonnegative constant $c$ that only depends on $\beta_{if}$.
\end{proposition}

\begin{proof}
Let $U_{f}=u_{f}^{1}-u_{f}^{2}$ be the difference of two solutions of
{\rm (\ref{Suf})} corresponding to $u_{s}^{1}$ and $u_{s}^{2}$ of $W_{z}\left(
T\right)  $, with same initial conditions at $z=0$ and set $U_{s}=u_{s}^{1}-u_{s}^{2}$ and $W_{f}\left(  r,z,t\right)
=U_{f}\left(  r,z,t\right)  -U_{s}\left(  z,t\right)  $ to obtain
\begin{equation}
\left\{
\begin{array}
[c]{rll}%
\dfrac{\partial W_{f}}{\partial z}-\beta_{f}\dfrac{1}{r\left(  1-r^{2}\right)
}\dfrac{\partial}{\partial r}\left(  r\dfrac{\partial W_{f}}{\partial
r}\right)  & = & -\dfrac{\partial U_{s}}{\partial z},\\
W_{f}\left(  r,0,t\right)  & = & 0,\\
W_{f}\left(  1,z,t\right)  & = & 0,\\
\dfrac{\partial W_{f}}{\partial r}\left(  0,z,t\right)  & = & 0.
\end{array}
\right.  \label{syst}%
\end{equation}
Multiplying {\rm (\ref{syst})} by $r(1-r^{2})W_{f}$, integrating on $\left]
0,1\right[  \times\left]  0,s\right[  \times\left]  0,\tau\right[  $ with
$s\in\left]  0,1\right[  $ and $\tau\in$ $\left]  0,T\right[  $ we obtain
using Young's inequality for $\varepsilon>0$ and $\int_{0}^{1}r\left(
1-r^{2}\right)  dr=1/4$ that 
\[%
\begin{array}
[c]{l}%
{\displaystyle\int\nolimits_{0}^{\tau}}
{\displaystyle\int\nolimits_{0}^{1}}
\left\Vert W_{f}\right\Vert ^{2}\left(  r,s,t\right)  r(1-r^{2})drdt+2%
{\displaystyle\int\nolimits_{0}^{\tau}}
{\displaystyle\int\nolimits_{0}^{s}}
{\displaystyle\int\nolimits_{0}^{1}}
\left\Vert \beta_{f}^{1/2}\dfrac{\partial W_{f}}{\partial r}\right\Vert
^{2}rdrdzdt\\
\leq\dfrac{\varepsilon}{4}%
{\displaystyle\int\nolimits_{0}^{\tau}}
{\displaystyle\int\nolimits_{0}^{1}}
\left\Vert \dfrac{\partial U_{s}}{\partial z}\right\Vert ^{2}dzdt+\dfrac
{1}{\varepsilon}%
{\displaystyle\int\nolimits_{0}^{\tau}}
{\displaystyle\int\nolimits_{0}^{s}}
{\displaystyle\int\nolimits_{0}^{1}}
\left\Vert W_{f}\right\Vert ^{2}r(1-r^{2})rdrdzdt
\end{array}
\]
and we deduce the following inequalities from Gronwall's lemma
\[%
\begin{array}
[c]{rll}%
{\displaystyle\int\nolimits_{0}^{\tau}}
{\displaystyle\int\nolimits_{0}^{s}}
{\displaystyle\int\nolimits_{0}^{1}}
\left\Vert W_{f}\right\Vert ^{2}r\left(  1-r^{2}\right)  drdzdt & \leq &
\dfrac{\varepsilon^{2}}{4}\left(  e^{s/\varepsilon}-1\right)
{\displaystyle\int\nolimits_{0}^{\tau}}
{\displaystyle\int\nolimits_{0}^{1}}
\left\Vert \dfrac{\partial U_{s}}{\partial z}\right\Vert ^{2}dzdt,\\%
{\displaystyle\int\nolimits_{0}^{\tau}}
{\displaystyle\int\nolimits_{0}^{s}}
{\displaystyle\int\nolimits_{0}^{1}}
\left\Vert \beta_{f}^{1/2}\dfrac{\partial W_{f}}{\partial r}\right\Vert
^{2}rdrdzdt & \leq & \dfrac{\varepsilon}{8}e^{s/\varepsilon}%
{\displaystyle\int\nolimits_{0}^{\tau}}
{\displaystyle\int\nolimits_{0}^{1}}
\left\Vert \dfrac{\partial U_{s}}{\partial z}\right\Vert ^{2}dzdt.
\end{array}
\]
Taking $\varepsilon=s$ we notice that $\frac{\partial W_{f}}{\partial r}%
=\frac{\partial U_{f}}{\partial r}$ because $U_{s}$ doesn't depend on $r$ to
have
\[%
{\displaystyle\int\nolimits_{0}^{T}}
{\displaystyle\int\nolimits_{0}^{1}}
{\displaystyle\int\nolimits_{0}^{1}}
\left\Vert \dfrac{\partial\left(  \Psi\left(  u_{s}^{1}\right)  -\Psi\left(
u_{s}^{2}\right)  \right)  }{\partial r}\right\Vert ^{2}rdrdzdt\leq\dfrac
{e}{8\inf_{i}\left(  \beta_{if}\right)  }%
{\displaystyle\int\nolimits_{0}^{T}}
{\displaystyle\int\nolimits_{0}^{1}}
\left\Vert \dfrac{\partial\left(  u_{s}^{1}-u_{s}^{2}\right)  }{\partial
z}\right\Vert ^{2}dzdt.
\]
From the fact that $U_{f}=W_{f}+U_{s}$ using previous inequalities we show
that
\[%
\begin{array}
[c]{l}%
\left\Vert U_{f}\right\Vert _{L^{2}\left(  \left]  0,T\right[  \times\left]
0,s\right[  ;L_{r(1-r^{2})}^{2}\left(  0,1\right)  \right)  ^{N}}\\

\begin{array}
[c]{l}%
\leq\sqrt{\dfrac{e-1}{4}%
{\displaystyle\int\nolimits_{0}^{\tau}}
{\displaystyle\int\nolimits_{0}^{1}}
\left\Vert \dfrac{\partial U_{s}}{\partial z}\right\Vert ^{2}dzdt}+\dfrac
{1}{2}\left\Vert U_{s}\right\Vert _{L^{2}\left(  \left]  0,T\right[  ;\left(
L^{2}\left(  0,s\right)  \right)  ^{N}\right)  },
\end{array}
\end{array}
\]
which leads to the announced property.
\end{proof}
\begin{theorem}
For $T$ small enough, {\rm (\ref{SPp_reg})} admits a solution in
\[
W\left(  T\right)  =\left\{
\begin{array}
[c]{r}%
\left(  u,v\right)  \in W_{r}\left(  T\right)  \times W_{z}\left(  T\right)
\mid u\left(  1,z,t\right)  =v\left(  z,t\right)  ,\qquad\\
\forall\left(  z,t\right)  \in\left]  0,1\right[  \times\left]  0,T\right[
\end{array}
\right\}  .
\]
\end{theorem}
\vspace{-2mm}
\begin{proof}
Using
{\rm (\ref{phi})} with $\Psi\left(  u_{s}^{1}\right)  $ and $\Psi\left(  u_{s}%
^{2}\right)  $ implies
\[%
\begin{array}
[c]{l}%
\left\Vert \Phi\left(  \Psi\left(  u_{s}^{1}\right)  \right)  -\Phi\left(
\Psi\left(  u_{s}^{2}\right)  \right)  \right\Vert _{W_{z}(T)}^{2}\\%
\begin{array}
[c]{l}%
\leq T\left(  aT+b\right)
{\displaystyle\int\nolimits_{0}^{T}}
{\displaystyle\int\nolimits_{0}^{1}}
\left\Vert \dfrac{\beta_{f}}{r\left(  1-r^{2}\right)  }\dfrac{\partial
}{\partial r}\left(  r\dfrac{\partial\left(  \Psi\left(  u_{s}^{1}\right)
-\Psi\left(  u_{s}^{2}\right)  \right)  }{\partial r}\right)  \right\Vert
_{W_{r0}^{\prime}}^{2}dzdt
\end{array}
\end{array}
\]
because of {\rm (\ref{syst})} with $a$ and $b$ nonnegative. This leads to
\begin{equation}
\left\Vert \Phi\left(  \Psi\left(  u_{s}^{1}\right)  \right)  -\Phi\left(
\Psi\left(  u_{s}^{2}\right)  \right)  \right\Vert _{W_{z}(T)}^{2}\leq
T\left(  aT+b\right)  \left\Vert \Psi\left(  u_{s}^{1}\right)  -\Psi\left(
u_{s}^{2}\right)  \right\Vert _{W_{r}(T)}^{2}. \label{in2}%
\end{equation}
Using {\rm (\ref{fluid})} with {\rm (\ref{in2})} we show that if $T$ is small enough the
mapping $\Phi\circ\Psi\ $is strictly contractant from $W_{z}(T)$ to itself
which prove the existence of a weak solution of {\rm (\ref{SPp_reg})} in the
apropriate space.
\end{proof}

\subsection{Back to the initial problem}

Let
\vspace{-2mm}
$$\widetilde{W}_{z=}\left\{  u\in\left(  L^{2}\left(  0,1\right)  \right)
^{N}\mid\dfrac{\partial u_{N}}{\partial z}\in L^{2}\left(  0,1\right)
\right\}  ,\text{ }\widetilde{W}_{z}\left(  T\right)  =\left\{  u\in
L^{2}\left(  0,T;\widetilde{W}_{z}\right)  \right\},$$
$$\widetilde{W}\left(  T\right)  =\left\{  \left(  u,v\right)  \in W_{r}\left(
T\right)  \times\widetilde{W}_{z}\left(  T\right)  \mid u\left(  1,z,t\right)
=v\left(  z,t\right)  ,\text{ }\forall\left(  z,t\right)  \in\left[
0,1\right]  \times\left[  0,T\right]  \right\},$$
we have the
\begin{theorem}
Assume $\theta_{is}$ goes to $0$ for $i\in\{1,\ldots,N-1\}$. The solution $\left(
C_{if},C_{is}\right)  $ of {\rm (\ref{SPp_reg})} weakly converge in $\widetilde
{W}\left(  T\right)  $ to the solution $\left(  \widetilde{C}_{if},\widetilde{C}_{is}\right)  $ of 
{\rm (\ref{sysp})}.
\end{theorem}

\begin{proof}
Take $\varphi\in$ $W\left(  T\right)  ,$ using {\rm (\ref{SPp_reg})} we have
\[%
\begin{array}
[c]{l}%
{\displaystyle\int\nolimits_{0}^{T}}
{\displaystyle\int\nolimits_{0}^{1}}
C_{if}\left(  r,1,t\right)  \varphi\left(  r,1,t\right)  r\left(
1-r^{2}\right)  drdt-%
{\displaystyle\int\nolimits_{0}^{T}}
{\displaystyle\int\nolimits_{0}^{1}}
C_{i0}\left(  r\right)  \varphi\left(  r,0,t\right)  r\left(  1-r^{2}\right)
drdt\\
-%
{\displaystyle\int\nolimits_{0}^{T}}
{\displaystyle\int\nolimits_{0}^{1}}
{\displaystyle\int\nolimits_{0}^{1}}
C_{if}\dfrac{\partial\varphi}{\partial z}r\left(  1-r^{2}\right)
drdzdt+\beta_{if}%
{\displaystyle\int\nolimits_{0}^{T}}
{\displaystyle\int\nolimits_{0}^{1}}
{\displaystyle\int\nolimits_{0}^{1}}
\dfrac{\partial C_{if}}{\partial r}\dfrac{\partial\varphi}{\partial
r}rdrdzdt\\
+\dfrac{\beta_{if}}{\gamma_{is}}%
{\displaystyle\int\nolimits_{0}^{1}}
C_{is}\left(  z,T\right)  \varphi\left(  1,z,T\right)  dz-\dfrac{\beta_{if}%
}{\gamma_{is}}%
{\displaystyle\int\nolimits_{0}^{1}}
C_{is0}\left(  z\right)  \varphi\left(  1,z,0\right)  dz\\
-\dfrac{\beta_{if}}{\gamma_{is}}%
{\displaystyle\int\nolimits_{0}^{T}}
{\displaystyle\int\nolimits_{0}^{1}}
C_{is}\dfrac{\partial\varphi}{\partial t}\left(  1,z,t\right)  dzdt+\dfrac
{\beta_{if}\theta_{is}}{\gamma_{is}}%
{\displaystyle\int\nolimits_{0}^{T}}
{\displaystyle\int\nolimits_{0}^{1}}
\dfrac{\partial C_{is}}{\partial z}\dfrac{\partial\varphi}{\partial z}\left(
1,z,t\right)  dzdt\\
-\dfrac{\beta_{if}}{\gamma_{is}}\delta_{i}%
{\displaystyle\int\nolimits_{0}^{T}}
{\displaystyle\int\nolimits_{0}^{1}}
\mathbf{r}_{i}(C_{1s},\ldots,C_{Ns})\varphi\left(  1,z,t\right)  dzdt=0.
\end{array}
\]
All the terms are bounded independently of $\theta_{is}$ (cf. Lemma
{\rm \ref{estimthe})}. To pass to the limit in the non-linear term with $\mathbf{r}%
_{i}$, we use the fact that $C_{is}$ is bounded in $L^{2}\left(  0,T;\left(
L^{2}\left(  0,1\right)  \right)  ^{N}\right)  $ and the Theorem {\rm 2.1} of Krasnoselskii
{\rm [8, page 22]}. This allows us to let $\theta_{is}$ going to $0$ for
$i\in\{1,\ldots,N-1\}$ in the precedent variationnal formulation to obtain
\[%
\begin{array}
[c]{l}%
{\displaystyle\int\nolimits_{0}^{T}}
{\displaystyle\int\nolimits_{0}^{1}}
\widetilde{C}_{if}\left(  r,1,t\right)  \varphi\left(  r,1,t\right)  r\left(
1-r^{2}\right)  drdt-%
{\displaystyle\int\nolimits_{0}^{T}}
{\displaystyle\int\nolimits_{0}^{1}}
C_{i0}\left(  r\right)  \varphi\left(  r,0,t\right)  r\left(  1-r^{2}\right)
drdt\\
-%
{\displaystyle\int\nolimits_{0}^{T}}
{\displaystyle\int\nolimits_{0}^{1}}
{\displaystyle\int\nolimits_{0}^{1}}
\widetilde{C}_{if}\dfrac{\partial\varphi}{\partial z}r\left(  1-r^{2}\right)
drdzdt+\beta_{if}%
{\displaystyle\int\nolimits_{0}^{T}}
{\displaystyle\int\nolimits_{0}^{1}}
{\displaystyle\int\nolimits_{0}^{1}}
\dfrac{\partial\widetilde{C}_{if}}{\partial r}\dfrac{\partial\varphi}{\partial
r}rdrdzdt\\
+\dfrac{\beta_{if}}{\gamma_{is}}%
{\displaystyle\int\nolimits_{0}^{1}}
\widetilde{C}_{is}\left(  z,T\right)  \varphi\left(  1,z,T\right)
dz-\dfrac{\beta_{if}}{\gamma_{is}}%
{\displaystyle\int\nolimits_{0}^{1}}
C_{is0}\left(  z\right)  \varphi\left(  1,z,0\right)  dz\\
-\dfrac{\beta_{if}}{\gamma_{is}}%
{\displaystyle\int\nolimits_{0}^{T}}
{\displaystyle\int\nolimits_{0}^{1}}
\widetilde{C}_{is}\dfrac{\partial\varphi}{\partial t}\left(  1,z,t\right)
dzdt+\delta_{N}^{i}\dfrac{\beta_{if}\theta_{is}}{\gamma_{is}}%
{\displaystyle\int\nolimits_{0}^{T}}
{\displaystyle\int\nolimits_{0}^{1}}
\dfrac{\partial\widetilde{C}_{is}}{\partial z}\dfrac{\partial\varphi}{\partial
z}\left(  1,z,t\right)  dzdt\\
-\dfrac{\beta_{if}}{\gamma_{is}}\delta_{i}%
{\displaystyle\int\nolimits_{0}^{T}}
{\displaystyle\int\nolimits_{0}^{1}}
\mathbf{r}_{i}(\widetilde{C}_{1s},\ldots,\widetilde{C}_{Ns})\varphi\left(
1,z,t\right)  dzdt=0.
\end{array}
\]
Making some particular choices for the test-function $\varphi$ 
we deduce that $\left(  \widetilde{C}_{if},\widetilde{C}_{is}\right)  $ is
solution of {\rm (\ref{SRp})} with the initial or boundary conditions {\rm (\ref{cond})}.
From {\rm (\ref{ineg_theta_is})} and the fact that the space of existence of
$\widetilde{C}_{if}$ doesn't depend on $\theta_{is}$ we have that the
solution $\left(  \widetilde{C}_{if},\widetilde{C}_{is}\right)  $ of
{\rm (\ref{SRp})} is in $\widetilde{W}\left(  T\right)$ and the existence of a solution of
{\rm (\ref{sysp})}.
\end{proof}

\begin{remark}
In the following we will write $\left(  C_{if},C_{is}\right)  $ instead of
$\left(  \widetilde{C}_{if},\widetilde{C}_{is}\right)  $ for the solution of
{\rm (\ref{sysp})}.
\end{remark}

\section{Uniqueness of the solution}
\vspace{-2mm}

\begin{proposition}
The system {\rm (\ref{sysp})} admits an unique solution $\left(  C_{if}%
,C_{is}\right)  _{i=1...N}$ in $\widetilde{W}\left(  T\right)  $.
\end{proposition}

\begin{proof}
Suppose there exists two solutions $\left(  C_{if}^{1},C_{is}^{1}\right)  $,
$\left(  C_{if}^{2},C_{is}^{2}\right)  $ of problem {\rm (\ref{sysp})}, and set 
$W_{if}=C_{if}^{1}-C_{if}^{2},$ $W_{is}=C_{is}^{1}-C_{is}^{2}.$
For $i\in\{1,...,N\}$, $\left(  W_{if},W_{is}\right)  $ is a weak solution of
\begin{equation}
\left\{
\begin{array}
[c]{rl}%
\dfrac{\partial W_{if}}{\partial z}\left(  r,z,t\right)  = & \beta_{if}%
\dfrac{1}{r\left(  1-r^{2}\right)  }\dfrac{\partial}{\partial r}\left(
r\dfrac{\partial W_{if}}{\partial r}\right)  \left(  r,z,t\right)  ,\\
\dfrac{\partial W_{is}}{\partial t}\left(  z,t\right)  = & -\gamma_{is}%
\dfrac{\partial W_{if}}{\partial r}(1,z,t)+\delta_{N}^{i}\theta_{Ns}%
\dfrac{\partial^{2}W_{Ns}}{\partial z^{2}}\left(  z,t\right) \\
& +\delta_{i}\left(  \mathbf{r}_{i}\left(  C_{1s}^{1},\ldots,C_{Ns}%
^{1}\right)  -\mathbf{r}_{i}\left(  C_{1s}^{2},\ldots,C_{Ns}^{2}\right)
\right)  \left(  z,t\right)  ,
\end{array}
\right.  \label{sys}%
\end{equation}
with the initial or boundary conditions
\[
\left\{
\begin{array}
[c]{rlrlrl}%
W_{if}\left(  r,0,t\right)  = & 0, &W_{if}\left(  1,z,t\right)  = & W_{is}\left(
z,t\right), & \dfrac{\partial W_{if}}{\partial
r}\left(  0,z,t\right)  = & 0 ,\\
\theta_{Ns}\dfrac{\partial W_{Ns}}{\partial
z}\left(  0,t\right)  = & 0, & \theta_{Ns}\dfrac{\partial W_{Ns}}{\partial
z}\left(  1,t\right)  = & 0, & W_{is}\left(  z,0\right)  = & 0.
\end{array}
\right.
\]
Multiplying {\rm (\ref{sys})$_\text{1}$} by $r\left(  1-r^{2}\right)  W_{if},$
integrating for $\tau\in\left]  0,T\right[  $, using {\rm (\ref{sys})$_\text{2}$}
and the fact that
the $\mathbf{r}_{i}$ are Lipschitz continuous with constant $k$ we
have\vspace{-1mm}
\[%
{\displaystyle\sum\limits_{i=1}^{N}}
{\displaystyle\int\nolimits_{0}^{1}}
\left(  W_{is}\right)  ^{2}\left(  z,\tau\right)  dz\leq2kN%
{\displaystyle\sum\limits_{i=1}^{N}}
{\displaystyle\int\nolimits_{0}^{\tau}}
{\displaystyle\int\nolimits_{0}^{1}}
\left(  W_{is}\right)  ^{2}dzdt.
\]
But $W_{is}(z,0)=0$ so we use Gronwall's lemma to conclude that $W_{is}%
(z,\tau)=0$ which implies that all the $W_{if}$ and $W_{is}$ are equal to $0$
in their existence spaces.
\end{proof}

\section{Conclusion}

Starting from a non-stationary model of catalytic converter with cylindrical passage presenting many 
mathematical difficulties we have used some parabolic regularization technics and constructed an appropriate 
functionnal space 
to prove that this problem 
admits one and only one solution for small time.

\vspace{-1mm}

\section*{Acknowledgment.} 
I would like to thank Professor Alain Brillard for giving me the subject of this paper for 
my PhD thesis, Professor Paul Deuring for his help in the proof of lemma 4, Professor Samir Akesbi for his advices and Professor Bernard Brighi for many useful discussions.

\end{document}